\documentclass[12pt]{amsart}
\usepackage{amsfonts, amsmath, amssymb, amscd, latexsym,graphicx}

\newtheorem{thm}{Theorem}[section]
\newtheorem{prop}[thm]{Proposition}
\newtheorem{lem}[thm]{Lemma}

\newtheorem{defn}[thm]{Definition}

\def\qed{\hfill\rlap{$\sqcup$}$\sqcap$\par}

\begin{document}

\title{Pure braid subgroups of braided Thompson's groups}

\author{Tom Brady}
\address{Department of Mathematics, Dublin City University, Dublin 9, Ireland } \email{Thomas.Brady@dcu.ie}
\thanks{The first author gratefully acknowledges the hospitality of the Centre de Recerca Matem\`atica.}

\author{Jos\'e Burillo}
\address{
Departament de Matem\'atica Aplicada IV, Universitat Polit\'ecnica
de Catalunya, Escola Polit{\`e}cnica Superior de Castelldefels,
08860 Castelldefels, Barcelona, Spain} \email{burillo@mat.upc.es}
\thanks{The second author acknowledges support from
NSF grant DMS-0305545 and the hospitality of the Centre de Recerca
Matem\`atica.}

\author{Sean Cleary}
\address{Department of Mathematics,
The City College of New York \& The CUNY Graduate Center, New
York, NY 10031} \email{cleary@sci.ccny.cuny.edu}\thanks{The third
author acknowledges support from PSC-CUNY grant \#66490, NSF grant
DMS-0305545 and the hospitality of the Centre de Recerca
Matem\`atica.}

\author{Melanie Stein}
\address{Department of Mathematics, Trinity College, Hartford, CT 06106} \email{melanie.stein@trincoll.edu }
\thanks{The fourth author gratefully acknowledges the hospitality of the Centre de Recerca Matem\`atica.}


\date{\today}



\begin{abstract}
We describe pure braided versions of Thompson's group $F$.  These
groups, $BF$ and $\widehat{BF}$, are subgroups of the braided
versions of Thompson's group $V$, introduced by Brin and Dehornoy.
Unlike $V$, elements of $F$ are order-preserving self-maps of the
interval and we use pure braids together with elements of $F$ thus
preserving order.  We define these groups and give normal forms
for elements and describe infinite and finite presentations of
these groups.
\end{abstract}

\maketitle

\section{Introduction and definitions}
Thompson's groups $F$ and $V$ have been studied from many
perspectives. Both groups can be understood as groups of locally
orientation-preserving piecewise-linear maps of the unit interval.
In the case of $F$, these maps are homeomorphisms, and in the case
of $V$ the maps are right-continuous bijections. In both cases the
breakpoints and discontinuities are restricted to be dyadic
rational numbers, and the slopes, when defined, are powers of 2.
Both groups can also be understood by means of rooted binary tree
pair diagrams---order-preserving in the case of $F$. Cannon, Floyd
and Parry \cite{cfp} give an excellent introduction to these
groups and several approaches to understanding their properties.

A rooted binary tree is a finite tree where every node has valence
3 except the root, which has valence two, and the leaves, which
have valence one.  We draw such trees with the root on top and the
nodes descending from it to the leaves along the bottom.
 The two nodes immediately
below a node are its {\it children\/}. A node and its two children
form a {\it caret\/}. A caret whose two children are leaves is
called an {\it exposed caret}. We number the leaves of a rooted
binary tree with $n$ carets and $n+1$ leaves from $0$ to $n$ in
any order, although we frequently choose to number the leaves of a
rooted tree in order from left to right.

\begin{figure}
\includegraphics[width=3in]{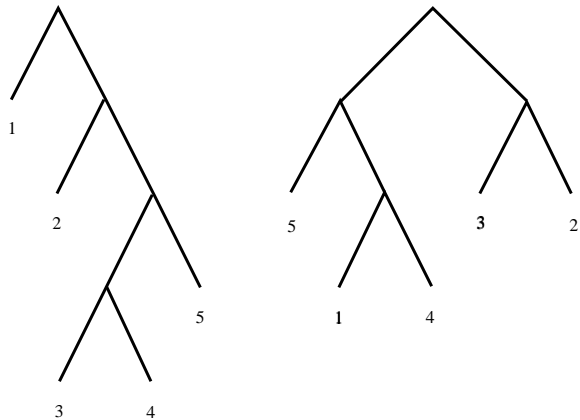}\\
\caption{An element of $V$\label{fig:velem}}
\end{figure}

A {\it tree pair diagram\/} is a triple $(T_-,\pi,T_+)$, where
$T_-$ and $T_+$ are two binary trees with the same number of
leaves $n$, and $\pi$ is a permutation in $S_n$. A reduction can
be performed in a diagram if both leaf numbers of an exposed caret
in $T_-$ are mapped by $\pi$ to the two leaf numbers of an exposed
caret in  $T_+$. In cases where a reduction is possible, we can
replace the exposed caret with a leaf and renumber leaves to give
an equivalent representative with a new permutation in $S_{n-1}$.
A tree pair diagram is {\it reduced} if no reductions are
possible. The set of binary tree pair diagrams thus admits an
equivalence relation, whose classes consist of those diagrams
which have a common reduced representative, with such reduced
representatives being unique. In Figure \ref{fig:velem} a reduced
diagram is depicted, where the leaf numbering describes
 the permutation.
The elements of $V$ which are actually in $F$ are precisely those
elements for which the permutation is the identity.

Composition in $V$ can be understood by means of these binary tree
diagrams. If two elements in $V$ are given by their representative
diagrams $(T_-,\pi,T_+)$ and $(S_-,\sigma,S_+)$, their composition
can be found by finding two tree diagrams in the corresponding
equivalence classes, $(T'_-,\pi',T'_+)$ and $(S'_-,\sigma',S'_+)$,
such that $T'_+=S'_-$. When this is achieved, the product element
is represented by the diagram $(T'_-,\pi'\circ\sigma',S'_+)$.

\begin{figure}
\includegraphics[width=1.5in]{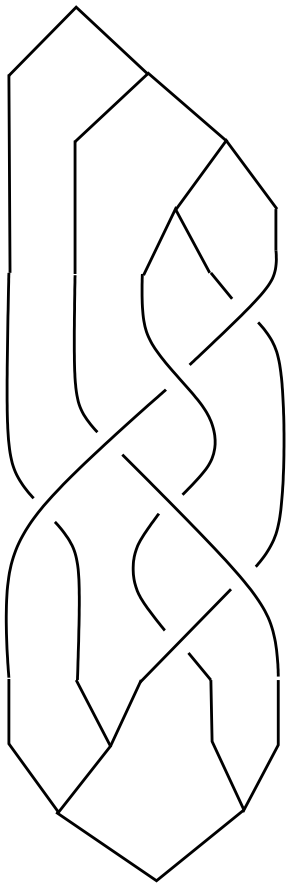}\\
\caption{An element of $BV$, one of the preimages of the element
drawn in Figure \ref{fig:velem} under the map
$\bar\phi$.\label{fig:bvelem}}
\end{figure}

Brin \cite{brinbv1, brinbv2}  and  Dehornoy \cite{dehornoy, dehornoy2}
describe braided Thompson's group $BV$, incorporating braids into
tree pair diagrams.
 The permutation  $\pi$ in a triple for an element of $V$ is replaced by a braid,
 giving the notion of a {\it braided tree diagram}. A braided tree diagram
is then a triple $(T_-,b,T_+)$, where the two trees have the same
number of leaves $n$, and $b\in B_n$ is a braid in $n$ strands.
Reductions can still be performed if two exposed carets are joined
by two parallel strands. In Figure \ref{fig:bvelem} an element of
$BV$ is depicted. We now draw the rooted tree $T_-$ with the root
at the top, and the tree $T_+$ below with the root at the bottom,
and then draw the braid between the leaves of the two trees as
indicated.

In \cite{brinbv1}, Brin also describes a ``larger" group
$\widehat{BV}$ as a group of braided forest diagrams having all
but finitely many of the forests trivial. Not only does
$\widehat{BV}$ naturally contain $BV$ as the subgroup of forest
diagrams where all trees are trivial except the first pair, but
$\widehat{BV}$ also sits inside $BV$ as the subgroup of braided
tree diagrams where the rightmost strand is always unbraided. We
begin this paper by providing finite and infinite presentations of
$BV$ which contain the presentations provided by Brin in
\cite{brinbv1} as subpresentations. Next, we describe subgroups
$BF$ of $BV$ and $\widehat{BF}$ of $\widehat{BV}$. Just as $F$ is
the subgroup of $V$ of order-preserving right continuous
bijections of $V$, the groups $BF$ and $\widehat{BF}$ are the
subgroups of order-preserving elements of the braided versions of
$V$.  The order is preserved by using generators which come from
$F$ and generators which involve pure braids.  We describe normal
forms for elements in these subgroups and obtain infinite and
finite presentations for these groups.  Dehornoy \cite{dehornoy}
calls this pure braid subgroup $PB_{\bullet}$, the group of pure
parenthesized braids.

\section{The braided Thompson's group $BV$}
In \cite{brinbv1}, an infinite presentation for $\widehat{BV}$ is
given. The generators in this presentation are the generators in
the standard infinite presentation for Thompson's group $F$, as
well as the generators for $B_{\infty}$. Here $B_{\infty}$ is
considered as a direct limit of the groups $B_n$, where $B_n$ is
included in $B_{n+1}$ via adding one strand at the right. Now
$\widehat{BV}$ sits naturally as a subgroup of $BV$; it is
isomorphic to the subgroup of all elements represented by braided
tree diagrams in which the rightmost strand is unbraided. Although
presentations for $BV$, both finite and infinite, are given in
\cite{brinbv2}, they are not related in a simple way to the
presentation for $\widehat{BV}$. Instead, we give a presentation
for $BV$ which contains  Brin's presentation for $\widehat{BV}$ as
a subpresentation. First, we define the set of generators. Any
element of $BV$ can be represented by a braided tree diagram
$(T_-, b, T_+)$ where both $T_-$ and $T_+$ have $n$ leaves and $b$
is a braid in $B_n$. A single tree can be thought of as a {\it
positive\/} element of Thompson's group $F$, when we take it as
being paired with an {\it all-right\/} tree, which is a tree whose
carets are all right children of their parent carets. These
positive elements correspond to elements which are positive words
in $F$ with respect to the infinite generating set $\{x_0, x_1,
\ldots\}.$ The correspondence between tree pair diagrams and
normal forms with respect to the infinite generating set is given
by the process of exponents of leaves, as described by Cannon,
Floyd and Parry \cite{cfp} and Fordham \cite{blakegd}. All-right
trees have all leaf exponents zero, and thus the normal forms for
tree pair diagrams which involve one all-right tree will be purely
negative or purely positive. We will denote by $R_n$ the all-right
tree which has $n$ leaves.

We can factor an element  $(T_-,b,T_+)$ into three pieces, using
all-right trees of the appropriate number of leaves, in a manner
similar to that done for elements of Thompson's group $T$ by
Burillo, Cleary, Stein and Taback \cite{thompt}.
 The resulting three elements in this factorization are
$$
(T_-,id,R_n)\qquad(R_n,b,R_n)\qquad(R_n,id,T_+),
$$ and the product of the elements represented by these three diagrams yields the original group element.
In general, these three tree pair diagrams will not be reduced; in
order for each of them to have the same number of carets, we may
need to take unreduced representatives for as many as two of the
three terms. By enlarging trees in this manner, it is clear that
every element of $BV$ can be factored this way.

Hence, we can always think  of an element of $BV$ as if it were
composed of two elements of Thompson's group $F$, one positive and
one negative, and one braid. It makes sense then to consider, as a
set of generators of $BV$, the set of generators for $F$ and the
set of generators for the braid groups, interpreted as braided
tree pair diagrams between all-right trees.

\begin{figure}
\includegraphics[width=3in]{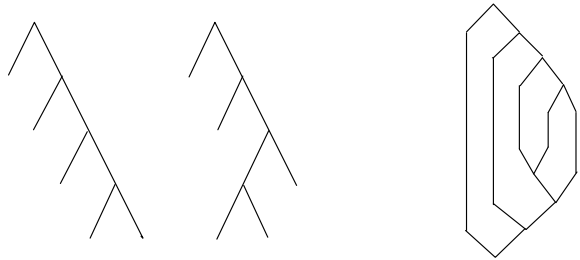}\\
\caption{The generator $x_2$ of $F$, in standard and in braided
form.\label{fig:x2}}
\end{figure}

The infinite set of generators for $F$ consists of the elements
$x_i$ with $i\ge 0$. Figure \ref{fig:x2} shows $x_2$ in both tree
pair diagram form and in braided tree form. These generators from
$F$ are enough to produce the two elements $(T_-,id,R_n)$ and
$(R_n,id,T_+)$ in $BV$.

We can consider the element $(R_n, b, R_n)$ as an element in the
appropriate braid group $B_n$.  Now this copy of $B_n$ is
generated by $n-1$ transpositions, the $i^{th}$ of which braids
strand $i$ over strand $i+1$.  We do not, however, need to include
all of these as generators. The generators which do not involve
braiding the last two strands can be obtained from the generators
of $B_{n-1}$ by splitting the last strand.

For this reason we must consider  two sets of generating braids,
one which leaves the rightmost strand unbraided and one which does
not. We define $\sigma_i$ to be the element represented by the
braided tree diagram $(R_{i+2}, a_i, R_{i+2})$, where $a_i$ is the
braid on $i+2$ strands which crosses strand $i$ over strand $i+1$.
Similarly, $\tau_i$ is the element represented by the diagram
$(R_{i+1}, b_i, R_{i+1})$, where $b_i$ is the braid on $i+1$
strands which crosses strand $i$ over stand $i+1$. Then the set
$\{\sigma_1, \sigma_2, \ldots, \sigma_{n-2},\tau_{n-1}\}$
generates the copy of $B_n$ containing all elements of $BV$
represented by diagrams of the form $(R_n, p, R_n)$. Notice that
the $x_i$ together with the $\sigma_i$ generate the copy of
$\widehat{BV}$ inside $BV$, and they correspond to Brin's
generators. We have shown:
\begin{prop}\label{bvgens} The elements $x_i$, for $i \geq 0$,
$\sigma_i$ for $i \geq 1$, and $\tau_i$ for $i \geq 1$ form a set
of generators for $BV$.
\end{prop}
There are three types of natural relations among these generators.
First, there are the generators involving only the generators of
$F$, namely $x_jx_i=x_ix_{i+j}$ for $j > i$. These are the
generators for the standard presentation for $F$.  Next, we expect
to need generators for each copy of $B_n$.  These yield four types
of relations:
\begin{itemize}
\item $ \sigma_i \sigma_j=\sigma_j \sigma_i$, for $j\geq i+2$
\item $\sigma_i \sigma_{i+1} \sigma_i = \sigma_{i+1} \sigma_i
\sigma_{i+1}$ \item $\sigma_i \tau_j = \tau_j \sigma_i$, for $j \geq
i+2$ \item $\sigma_i \tau_{i+1} \sigma_i = \tau_{i+1} \sigma_i
\tau_{i+1}$.
\end{itemize}
Finally, there are relations governing the interactions between
the generators for $F$ and the generators for the braid groups
$B_n$.
\begin{itemize}
\item $\sigma_i x_j=x_j\sigma_i$, for $i<j$ \item $\sigma_i
x_i=x_{i-1} \sigma_{i+1} \sigma_{i}$ \item $\sigma_i x_j=x_j
\sigma_{i+1}$,for $i \geq j+2$ \item
$\sigma_{i+1}x_i=x_{i+1}\sigma_{i+1}\sigma_{i+2}$ \item $\tau_i x_j
= x_j \tau_{i+1}$, for $i\geq j+2$ \item $\tau_i x_{i-1}= \sigma_i
\tau_{i+1}$ \item $\tau_i=x_{i-1}\tau_{i+1}\sigma_i$.
\end{itemize}

In preparation for showing that the relations above give a
presentation, we first introduce a special class of words in the
generators. We would like to identify those words in the
generators which could be identified easily with a triple of
diagrams in $BV$. As noted earlier, any element of $BV$ can be
represented by a triple of braided tree diagrams of the form
$(T_-, id, R_n)~~~(R_n,b,R_n)~~~(R_n, id, T_+)$. Such a triple
leads easily to a word in the generators as follows. The group
element represented by the first diagram is a positive element $a$
of $F \subset BV$, and may be expressed uniquely as a word of the
form $x_{i_1}^{r_1}x_{i_2}^{r_2} \ldots x_{i_k}^{r_k}$, where $i_1
< i_2 < \cdots < i_k$ and $r_m \geq 1$ for all $m$. Similarly, the
group element $c\in F \subset BV$ represented by the third diagram
can be uniquely expressed as a word of the form $x_{j_l}^{-s_l}
\ldots x_{j_2}^{-s_2} x_{j_1}^{-s_1}$, where $j_1 < j_2 < \cdots <
j_l$ and $s_m \geq 1$ for all $m$. Now the group element
represented by the middle diagram may be represented as some word
in the generators $\sigma_1, \sigma_2, \ldots,
\sigma_{n-2}$,$\tau_{n-1}$, and their inverses. For convenience,
we will call any word in this set of generators and their inverses
{\it{a word in the $B_n$ generators}}. Note that if such a word
contains no $\tau$ generators, it can be considered a word in the
$B_n$ generators for many values of $n$. Notice that the minimum
number of carets required in the trees for tree pair diagrams
representing $a$ and $c$ respectively, is at most $n-1$. The
concatenation of the three words described above yield a word
which cannot serve as a normal form, since we have not specified
preferred arrangements of the  $\sigma$'s and $\tau$ and
furthermore, there are many different triples of tree pair
diagrams representing any element. However, these words are nice
in that any word of the above special form can be easily
translated into a triple of diagrams, and we find them to be
useful tools.

Given a word $w\in F$, denote by $N(w)$ the number of carets in
the reduced binary tree diagram representing it. Here is the
algebraic description of blocks, which are these words which come
from a single triple of diagrams.

\begin{defn}
A word in the generators $x_i^{\pm 1}, \sigma_i^{\pm 1},
\tau_i^{\pm 1}$ is called a block if it is of the form
$w_1w_2w_3^{-1}$ where \begin{itemize} \item [(1)] $w_1$ is of the
form $x_{i_1}^{r_1}x_{i_2}^{r_2} \ldots x_{i_k}^{r_k}$, where $i_1
< i_2 < \cdots < i_k$ and $r_m \geq 1$ for all $m$. \item [(2)]
$w_3$ is of the form $x_{j_1}^{s_1} x_{j_2}^{s_2} \ldots
x_{j_l}^{s_l}$ , where $j_1 < j_2 < \cdots < j_l$ and $s_m \geq 1$
for all $m$, and by $w_3^{-1}$ we mean the word $x_{j_l}^{-s_l}
\ldots x_{j_2}^{-s_2}  x_{j_1}^{-s_1}$. \item [(3)] Let
$N=max(N(w_1), N(w_2))$. Then there exists an integer $n$, $n \geq
N+1$, such that $w_2$ is a word in the $B_n$ generators.
\end{itemize}
\end{defn}
Then we have the following lemma:
\begin{lem}\label{blockid}
A block $w_1w_2w_3^{-1}$ is the identity in $BV$ if and only if
$w_1$ and $w_3$ are the same word, and $w_2$ is the identity in
the copy of the braid group $B_n$ generated by $\tau_{n-1}^{\pm
1}$ and $\sigma_i^{\pm 1}$ where $1 \leq i \leq n-2$. \end{lem}
\begin{proof}
The lemma follows directly from the fact that any word which is a
block $w_1w_2w_3^{-1}$ can be represented by a braided tree
diagram $(T_-, b, T_+)$ where $w_1$ is represented by $(T_-, id,
R_n)$, $w_2$ is represented by $(R_n, b, R_n)$, and $w_3^{-1}$ is
represented by $(R_n, id, T_+)$, and from any such triple of
diagrams a block can be read off, unique up to the choice of the
word in $B_n$ expressing $b$. Since the identity in $BV$ can be
represented by the diagram consisting of the tree with only one
vertex, and the trivial braid on one strand, all other diagrams
representing the identity result from splitting strands, and will
always have two identical trees with the trivial braid. But such
diagrams translate into blocks of the form described in the lemma.
\end{proof}
We will use these blocks to prove:
\begin{thm}\label{bvinfpres}
The group $BV$ admits a presentation with generators:
\begin{itemize}
\item $x_i$, for $i \geq 0$, \item $\sigma_i$, for $i \geq 1$,
\item $\tau_i$, for $i \geq 1$.
\end{itemize}
and relators
\begin{itemize}
\item [(A)]$x_j x_i=x_i x_{j+1} \text{  for }j > i$ \item
[(B1)]$\sigma_i \sigma_j=\sigma_j \sigma_i \text{  for }j-i \geq
2$ \item [(B2)]$\sigma_i \sigma_{i+1} \sigma_i = \sigma_{i+1}
\sigma_i \sigma_{i+1}$ \item [(B3)]$\sigma_i \tau_j = \tau_j
\sigma_i\text{ for }j-i \geq 2$ \item [(B4)]$\sigma_i \tau_{i+1}
\sigma_i = \tau_{i+1} \sigma_i \tau_{i+1}$. \item [(C1)]$\sigma_i
x_j=x_j\sigma_i,\text{ for  }i<j$ \item [(C2)]$\sigma_i
x_i=x_{i-1} \sigma_{i+1} \sigma_{i}$ \item [(C3)]$\sigma_i x_j=x_j
\sigma_{i+1},\text{ for  }i \geq j+2$ \item
[(C4)]$\sigma_{i+1}x_i=x_{i+1}\sigma_{i+1}\sigma_{i+2}$ \item
[(D1)]$\tau_i x_j = x_j \tau_{i+1},\text{ for  }i-j \geq 2$ \item
[(D2)]$\tau_i x_{i-1}= \sigma_i \tau_{i+1}$ \item
[(D3)]$\tau_i=x_{i-1}\tau_{i+1}\sigma_i$.
\end{itemize}
\end{thm}

This presentation appears without proof in J. Belk's thesis
\cite{belkdiss}.

\begin{proof}
Let $G$ be the abstract group given by the presentation above. We
map $G$ to $BV$ via $\phi$ by sending each generator to the
element of $BV$ with the same name. All relations in the
presentation hold in $BV$, so $\phi$ is a well-defined
homomorphism. Proposition \ref{bvgens} shows that the map is
surjective, so it remains only to show that $\phi$ is injective.
To show this, we must show that any word in the generators which
maps to the identity in $BV$ is already the identity in $G$. Now
by Lemma  \ref{blockid}, this is true if the word in question
happens to be a block. So we are done once we show that the
relations in $G$ are sufficient to transform any word into a
block. But since any generator is itself a block, an arbitrary
word of length $k$ is trivially the product of $k$ blocks. So to
show $\phi$ is injective it is sufficient to prove that a word in
$G$ which is the product of two blocks can be rewritten, using the
relations in $G$, as a single block. We first prove a series of 3
preliminary lemmas, from which we will deduce this fact in Lemma
\ref{2blockstoone}, which will complete the proof of the theorem.
\end{proof}
Our first lemma permits us to push the $x_i^{\pm 1}$ generators to
the left or right of the braid generators, which helps move a word
toward block form.
\begin{lem}\label{push}
If $w$ is a word in the $B_n$ generators, and $i \leq n-2$, then
$x_i^{-1} w $ is equivalent in $G$ to either $\bar{w} x_{i'}^{-1}$
or $\bar{w}$, where $\bar{w}$ is a word in the $B_{n+1}$
generators, and $i' \leq n-2$. Similarly, under the same
conditions on all indices, a word $wx_i$ may be replaced by either
$x_{i'}\bar{w}$ or $\bar{w}$.
\end{lem}

\begin{proof}
We describe first how to push $x_i^{-1}$ past the $\sigma$ and
$\tau$ generators. Pushing past $\sigma$ type generators is always
possible, but in order to push past $\tau$'s we must carefully
keep track of the index of the $x^{-1}$ as it moves along. Using
the relations of type C, we may replace $x_k^{-1}\sigma_j^{\pm1}$
by $w(\sigma)x_{k'}^{-1}$, where $w(\sigma)$ is a word in the
$\sigma$ generators and their inverses of length one or two.
Furthermore, the maximum index appearing in $w(\sigma)$ is $j+1$.
Now the index $k'$ can, in general, be either $k, k-1$, or $k+1$.
However, it only increases to $k+1$ in the case where $j=k+1$
also. So since the initial index $i$ satisfies $i \leq n-2$ and $j
\leq n-2$, even a series of such replacements results in the
presence of $x_{i'}^{-1}$ with $i' \leq n-2$. This is important,
since relations $(D2)$ and $(D3)$ allow replacement of
$x_{n-2}^{-1} \tau_{n-1}^{\pm1}$ by
$\tau_n^{\pm1}\sigma_{n-1}^{\pm1}$, and relations $(D1)$ allow us
to replace $x_k^{-1}\tau_{n-1}^{\pm1}$ by $\tau_n^{\pm1}x_k^{-1}$
if $k \leq n-3$.  Hence, $x_i^{-1}w$ can be replaced by either
$\bar{w}x_{i'}^{-1}$ or simply $\bar{w}$ as claimed.  The argument
for $wx_i$ is similar.
\end{proof}
Next we prove a lemma showing that a word in the $B_n$ generators
can always be pumped up to a word in the $B_{n+1}$ generators at
the expense of tacking on an $x$ generator.
\begin{lem}\label{pumping}
Let $w$ be a word in the $B_n$ generators. Then using the relators
in $G$, $w$ may be replaced by either $\bar{w}$ or
$\bar{w}x_i^{-1}$ where $i \leq n-2$ and $\bar{w}$ is a word in
the $B_{n+1}$ generators. Similarly, $w$ may also be replaced by
either $\bar{w}$ or $x_i \bar{w}$ where $\bar{w}$ is a word in the
$B_{n+1}$ generators.
\end{lem}
\begin{proof} Consider the leftmost occurrence of $\tau_{n-1}$ in the
word $w$, that is, $w=w_1 \tau_{n-1}^{\pm1}w_2$, where $w_1$ has
only $\sigma$ generators. Using relations $(D2)$ or $(D3)$
depending on the exponent of $\tau_{n-1}$, replace $w$ by $w_1
\sigma_{n-1}^{\pm1} \tau_n^{\pm1}x_{n-2}^{-1}w_2$, and then apply
Lemma \ref{push} to $x_{n-2}^{-1}w_2$ to replace it with either
$\bar{w}_2x_i^{-1}$ or $\bar{w}_2$ with $i \leq n-2$ and where
$\bar{w}_2$ is a word in the $B_{n+1}$ generators. Then the
desired $\bar{w}$ is
$w_1\sigma_{n-1}^{\pm1}\tau_n^{\pm1}\bar{w}_2$. Similarly, working
from the right, $w$ can be replaced by either $x_i\bar{w}$ or
$\bar{w}$.
\end{proof}
 The two previous lemmas will now be used to show that the relators allow us to transform the product of two blocks to a new product of two blocks where the combined length of the middle two of the 6 subwords involved is reduced.
\begin{lem}\label{squish}
Let $w=w_1w_2w_3^{-1}$ and $v=v_1v_2v_3^{-1}$ be two blocks. Then
the relations in $G$ allow us to replace the word $wv$ by
$w_1'w_2'w_3'^{-1}v_1'v_2'v_3'^{-1}$, the product of two blocks
$w'=w_1'w_2'w_3'^{-1}$ and $v'=v_1'v_2'v_3'^{-1}$, where
$l(v_1')+l(w_3') < l(v_1)+l(w_3)$.
\end{lem}
\begin{proof}
Let $x_{i_w}$ and $x_{i_v}$ be the first letters in $w_3$ and
$v_1$. If they are the same, we can delete the pair
$x_{i_w}^{-1}x_{i_v}$ and we are done. If not, suppose $i_w < i_v$
(if $i_v < i_w$ a similar argument works, truncating $v_3$ and
absorbing $x_{i_v}$ into $w$). Let $w_1'=w_1$, $w_2'=w_2$, and let
$w_3'$ be $w_3$ with $x_{i_w}$ deleted. Now we use the relations
(A) to replace $x_{i_w}^{-1}v_1$ by $v_1'x_{i_w}^{-1}$. Note that
$v_1'$ and $v_1$ have the same length, but $N(v_1')=N(v_1)+1$,
since each index in $v_1$ is increased by $1$ as $x_{i_w}^{-1}$
moves past it (see Theorem 3 of \cite{bcs}). Next, suppose $v_2$
is a word in $\tau_{n-1}^{\pm1}$ and $\sigma_j^{\pm1}$ with $1
\leq j \leq n-2$. Then $N(v_1) \leq n$ and $N(v_3) \leq n$, so
$N(v_1') \leq n+1$.  We must replace $x_{i_w}^{-1}v_2v_3^{-1}$ by
$v_2'v_3'^{-1}$ so that $v_1'v_2'v_3'^{-1}$ is a block. We will
consider two cases.

{\it Case 1:} If $i_w \leq n-2$, we use Lemma \ref{push} to
replace $x_{i_w}^{-1}v_2$ by $v_2'x_i^{-1}$ with $i \leq n-2$ and
$v_2'$ a word in $\tau_n^{\pm1}$ and $\sigma_j^{\pm1}$. Then
$(v_3x_i)^{-1}$ can be rewritten using relations (A) as a word
$(v_3')^{-1}$ so that $v_3'$ is a positive word in the generators
of $F$ with increasing indices from left to right. Then it again
follows from \cite{bcs} that $N(v_3')=Max(N(v_3)+1, i+2)$. Hence
$N(v_3') \leq n+1$, since $N(v_3)+1 \leq n+1$ and $i+2 \leq n$,
and this implies that $v_1'v_2'(v_3')^{-1}$ is a block as desired.

{\it Case 2:} If $i_w > n-2$, it is necessary to first use Lemma
\ref{pumping} to replace $v_2$ by either $\bar{v}_2$ or
$\bar{v}_2x_i^{-1}$ where $i \leq n-2$. If $x_i^{-1}$ is present,
we use relations (A) to replace $(v_3x_i)^{-1}$ by
$\bar{v}_3^{-1}$, and we see that $N(\bar{v}_3)=Max(N(v_3)+1, i+2)
\leq n+1$.  We continue applying Lemma \ref{pumping} and absorbing
any resulting $x_i^{-1}$ letters into the $v_3^{-1}$ part of the
word in this manner, and after $i_w-(n-2)$ repetitions we have
replaced $x_{i_w}^{-1}v_1v_2v_3^{-1}$ by
$v_1'x_{i_w}^{-1}\bar{v}_2\bar{v}_3^{-1}$, where $\bar{v}_2$ a
word in the $B_{i_w+2}$ generators, and $\bar{v}_3$ is a word in
the $x_i$ with indices increasing from left to right with
$N(\bar{v}_3) \leq n+(i_w-(n-2))=i_w+2$.  Now just as before we
can apply Lemma \ref{push} to replace
$v_1'x_{i_w}^{-1}\bar{v}_2\bar{v}_3^{-1}$ by
$v_1'v_2'x_i^{-1}\bar{v}_3^{-1}$ where $i \leq i_w$, and $v_2'$ is
a word in the $B_{i_w+3}$ generators. When we use relations (A) to
replace $(\bar{v}_3x_i)^{-1}$ by $v_3'^{-1}$, $N(v_3') =
Max(N(\bar{v}_3)+1, i+2)$. But since $N(\bar{v}_3)+1 \leq
n+(i_w-(n-2))=i_w+3$ and $i+2 \leq i_w+2$, $N(v_3') \leq i_w +3$,
and hence $v_1'v_2'v_3'$ is a block.

\end{proof}

Now we are in a position to prove the final lemma which completes
the proof of Theorem \ref{bvinfpres}.
\begin{lem}\label{2blockstoone}
The product of two blocks may be rewritten, using the relations of
$G$, as a single block.
\end{lem}
\begin{proof}
Let $w=w_1w_2w_3^{-1}$ and $v=v_1v_2v_3^{-1}$ be two blocks. We
apply Lemma \ref{squish}, at most $l(w_3)+l(v_1)$ times, to
replace $wv$ by $w_1'w_2'v_2'(v_3')^{-1}$ where $w_2'$ is a word
in the $B_n$ generators, $v_2'$ is a word in the $B_{n'}$
generators, $N(w_1') \leq n$, and $N(v_3') \leq n'$.  If $n=n'$,
declaring $u_1=w_1'$, $u_2=w_2'v_2'$, and $u_3=v_3'$ shows that
$u_1u_2u_3^{-1}=w_1'w_2'v_2'(v_3')^{-1}$ is a block.  If not, say
$n'<n$, we apply Lemma \ref{pumping} $n-n'$ times to replace
$v_2'$ by a word $\bar{v}_2$, a word in the $B_n$ generators,
followed by some new $x^{-1}$ generators, so that $v_2'v_3'^{-1}$
has been replaced by $\bar{v}_2(v_3'x_{i_1}x_{i_2} \cdots
x_{i_{n-n'}})^{-1}$ where $i_j \leq n'+j-3$ for $1 \leq j \leq
n-n'$. Now $N(v_3'x_{i_1})=Max(N(v_3)+1,i_1+2) \leq n'+1$. So
inductively, we have that $N(v_3'x_{i_1}x_{i_2} \cdots
x_{i_{n-n'}}) \leq n$, and hence declaring $u_1=w_1'$,
$u_2=w_2'\bar{v}_2$, and $u_3=\bar{v}_3$, shows that
$u_1u_2u_3^{-1}=w_1'w_2'\bar{v}_2\bar{v}_3^{-1}$ is a block. Of
course, if $n' > n$, a similar argument works, pumping up indices
in $w_2'$ instead.
\end{proof}
We remark that an infinite presentation of Thompson's group $V$ is
easily obtained from the presentation for $BV$ by adding two more
infinite families of relators, $\sigma_i^2=1$ and $\tau_i^2=1$ for
all $i$.

\section{A finite presentation for $BV$}
It is common for these type of infinite presentations for
Thompson-type groups to reduce to finite presentations. For
example, in \cite{cfp} an inductive argument is spelled out which
obtains the standard two generator-two relator presentation for
$F$ from the standard infinite presentation.  Brin uses similar
arguments to obtain finite presentations for $BV$ and
$\widehat{BV}$ from his infinite ones. In a similar manner, the
infinite presentation for $BV$ in the previous section reduces to
a finite presentation with 4 generators and only 18 relators, an
improvement over the presentation in \cite{brinbv2}, which has 4
generators and 26 relators.
\begin{thm} \label{finitebv}
The group $BV$ admits a finite presentation with generators $x_0,
x_1, \sigma_1, \tau_1$ and relators
\begin{itemize}\item[(a)] $x_2x_0=x_0x_3$, $x_3x_1=x_1x_4$
\item[(c1)] $\sigma_1x_2=x_2\sigma_1$, $\sigma_1x_3=x_3\sigma_1$,
$\sigma_2x_3=x_3\sigma_2$, $\sigma_2x_4=x_4 \sigma_2$ \item[(c3)]
$\sigma_2x_0=x_0\sigma_3$, $\sigma_3x_1=x_1 \sigma_4$ \item[(c4)]
$\sigma_1x_0=x_1\sigma_1\sigma_2$, $\sigma_2x_1=x_2\sigma_2\sigma_3$
\item[(d1)] $\tau_2x_0=x_0\tau_3$, $\tau_3x_1=x_1\tau_4$ \item[(d2)]
$\tau_1x_0=\sigma_1 \tau_2$, $\tau_2x_1=\sigma_2\tau_3$
\item[(b1)] $\sigma_1\sigma_3=\sigma_3\sigma_1$ \item[(b2)]
$\sigma_1\sigma_2\sigma_1=\sigma_2 \sigma_1 \sigma_2$ \item[(b3)]
$\sigma_1\tau_3=\tau_3\sigma_1$ \item[(b4)] $\sigma_1 \tau_2
\sigma_1=\tau_2 \sigma_1 \tau_2$
\end{itemize}
where the letters in the relators not in the set of 4 generators
are defined inductively by $x_{i+2}=x_i^{-1} x_{i+1}x_i$ for $i
\geq 0$, $\sigma_{i+1}=x_{i-1}^{-1}\sigma_ix_i\sigma_i^{-1}$ for
$i \geq 1$, and $\tau_{i+1}=x_{i-1}^{-1}\tau_i\sigma_i^{-1}$ for
$i \geq 1$.
\end{thm}
\begin{proof}
That the two (a) relators yield inductively all (A) relators in
the infinite presentation is a standard argument, given in
\cite{cfp}. Notice that the relators (C2) and (D3) in the infinite
presentation are precisely the relations used to inductively
define the higher index generators in the infinite presentation.
Now a straightforward induction yields the (C1) relators in the
infinite presentation from the (c1) relators, then the (C3)
relators from the (c3) relators, and so on, in the order the
groups of relators are listed in the finite presentation above. As
an example, we spell out the induction for the (B3) relators. So
suppose we have (b3), or $\sigma_1\tau_3=\tau_3\sigma_1$. Then
suppose inductively that we have established
$\sigma_1\tau_i=\tau_i\sigma_1$ for $3 \leq i <k$, where $k \geq
4$. Then $\sigma_1\tau_k=\sigma_1(x_{k-2}^{-1}\tau_{k-1}
\sigma_{k-1}^{-1})$ using the relator defining $\tau_k$. Now we
can move the $\sigma_1$ to the right, first using the (C1)
relators, then the inductive hypothesis, and finally the (B1)
relators, and then use the defining relation for $\tau_k$ in the
other direction, to obtain the relator
$\sigma_1\tau_k=\tau_k\sigma_1$. Therefore, by induction,
$\sigma_1\tau_k=\tau_k\sigma_1$ for all $k \geq 3$.  Now suppose
that we have $\sigma_i\tau_j=\tau_j\sigma_i$ for $j-i \geq 2$, $1
\leq i <k$, and $k \geq 2$. Then it follows that
$\sigma_k\tau_j=\tau_j\sigma_k$ for $j-k \geq 2$. We replace
$\sigma_k$ in the word $\sigma_k\tau_j$ by $x_{k-2}^{-1}
\sigma_{k-1} x_{k-1} \sigma_{k-1}^{-1}$, and then moves the
$\tau_j$ to the left, first using the inductive hypothesis to
obtain $x_{k-2}^{-1}\sigma_{k-1}x_{k-1}\tau_j\sigma_{k-1}^{-1}$,
then using the (D1) relators to obtain
$x_{k-2}^{-1}\sigma_{k-1}\tau_{j-1}x_{k-1}\sigma_{k-1}^{-1}$, and
finally using the inductive hypothesis again to obtain
$x_{k-2}^{-1}\tau_{j-1}\sigma_{k-1}x_{k-1}\sigma_{k-1}^{-1}$. Now
use (D1) relators to move $\tau_{j-1}$ left, to obtain
$\tau_jx_{k-2}^{-1}\sigma_{k-1}x_{k-1}\sigma_{k-1}^{-1}$. But now
the rightmost four letters can be replaced by $\sigma_k$ using the
defining relation for $\sigma_k$ in reverse, showing that
$\sigma_k\tau_j=\tau_j\sigma_k$. Hence, by induction, all (B3)
relators hold.

\end{proof}
For the corresponding finite presentation for $V$, we note that
the relations $\tau_i^2=1$ and $\sigma_i^2=1$ for $i \geq 2$ can
be deduced inductively from the two relations
$\tau_1^2=\sigma_1^2=1$, using the (C2) and (C3) relators in the
case of $\sigma_i$, and the (D2) and (D3) relators for $\tau_i$.
This yields a presentation for $V$ with 4 generators and 20
relators, not quite as efficient as the presentation in \cite{cfp}
with 14 relators.

\section{The group $PBV$}

The pure braid groups $P_n$ are the groups of braids where the
$i$th strand is braided with the other strands but returns to the
$i$th position.  There are several possible ways to construct
analogous subgroups of $BV$.    One way is by considering the
standard short exact sequences for the braid groups, involving the
pure braid groups and the permutation groups.  For each $n$, we
have:
$$
\CD 1 @>>> P_n @>>> B_n @>\phi_n>> S_n @>>> 1
\endCD
$$
which maps a braid to its permutation, and whose kernel is the
pure braid group $P_n$.  This family of maps $\phi_n$ collectively
induces a map
$$
\bar\phi: BV\longrightarrow V
$$
defined by $\bar\phi(T_-,b,T_+)=(T_-,\phi_n(b),T_+)$, where we use
the appropriate $\phi_n$ for the number of leaves in either tree.

Let $PBV=\ker \bar\phi$. By definition, a diagram $(T_-,b,T_+)$
represents an element in $PBV$ if it maps to the identity in $V$,
that is, if $\phi(b)=\text{id}$, and if $T_-=T_+$. Hence, $PBV$ is
the subgroup of $BV$ which consists of those elements which admit
a representative $(T,p,T)$ on which the two trees are the same and
the braid is pure. If an element admits one representative where
the two trees are equal, then every representative will have  the
trees being equal. So $PBV$ is  a subgroup, because the product of
such two elements also has representatives where the two trees are
equal. Note that it is crucial in this construction that the braid
is a pure braid.

The main result concerning the group $PBV$ is the following.

\begin{thm}\label{PBVnotfp} The group $PBV$ is not finitely generated.
\end{thm}

{\it Proof.\/} Given two elements of $PBV$ by their diagrams
$(T,p,T)$ and $(S,q,S)$, their product always admits a
representative diagram $(R,r,R)$, where $R$ is the {\it least
common multiple\/} of $S$ and $T$-- that is, the minimal tree
which contains both $S$ and $T$ as subtrees. Hence, if $PBV$ were
to be generated by a finite set $(T_i,p_i,T_i)$, for
$i=1,\ldots,k$, every tree in $PBV$ would admit a representative
whose tree would be the least common multiple of the $T_i$. There
are elements whose smallest representatives are of increasing size
and thus  $PBV$ cannot be finitely generated.\qed

\section{The braided Thompson's group $BF$}

From the map $\bar\phi$ defined above, since $F$ is the subgroup
of $V$ of those elements whose permutation is the identity, we can
define the group $BF=\bar\phi^{-1}(F)$. The group $BF$ is the
subgroup of $BV$ of those elements which admit a representative
$(T_-,p,T_+)$, where $p$ is a pure braid. Note the contrast with
$PBV$, since here the two trees are not necessarily equal. In
fact, $PBV$ is a subgroup of $BF$, and observing the restriction
to $BF$ of the map $\bar\phi$ above, it is easy to see that $PBV$
is also the kernel of $\bar\phi\big|BF$. Thus the diagram below is
commutative:

$$
\CD
1@>>>PBV@>>>BF@>\bar\phi|BF>>F@>>>1\\
@.    @|     @VVV  @VVV      \\
1@>>>PBV@>>>BV@>\bar\phi>>V@>>>1
\endCD
$$

The main goal of the remainder of this paper is to prove that $BF$
is finitely presented and to find both finite and infinite
presentations.

Finding generators for $BF$ is not difficult. Just as for $BV$, an
element of $BF$ is given by a triple $(T_-,p,T_+)$, where this
time $p$ is a pure braid. Again, we factor the element into three
pieces
$$
(T_-,id,R_n)\qquad(R_n,p,R_n)\qquad(R_n,id,T_+),
$$
where the individual diagrams may not be reduced. Hence, we can
always think of an element of $BF$ as if it were composed of two
elements of Thompson's group $F$, one positive and one negative,
and one pure braid. Again, just as for $BV$, we take as a set of
generators of $BF$, the set of generators for $F$ and the set of
generators for the pure braid groups, interpreted as braided tree
pair diagrams between all-right trees. We consider the element
$(R_n,p,R_n)$ as an element of the appropriate group $P_n$ of pure
braids. To generate these groups $P_n$ we would like to use the
braids $A_{ij}$, for $i< j$, which wrap the $i$-th strand around
the $j$-th one. See Hansen \cite{Hansen} for details of these
generating sets. The process of obtaining the generators of $BF$
from the generators $A_{ij}$ of $P_n$ is the same than the process
specified above for $BV$ from the standard braid generators.

\begin{figure}
\includegraphics[width=3in]{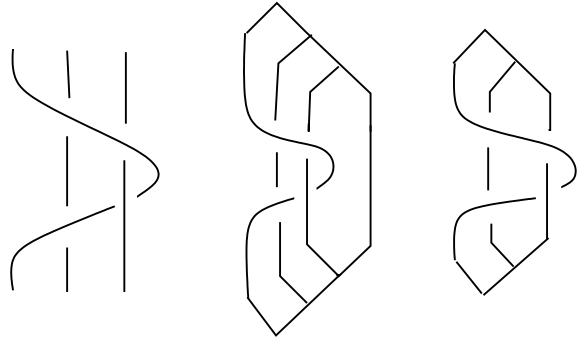}\\
\caption{The generator $A_{13}$ in $P_3$, with its corresponding
generators $\alpha_{13}$ and $\beta_{13}$ in $BF$.
\label{fig:a13}}
\end{figure}

We will denote by $\alpha_{ij}$ the element
$(R_{j+1},A_{ij},R_{j+1})$, and by $\beta_{ij}$ the element
$(R_j,A_{ij},R_j)$. As in $BV$, the differences between these two
families of generators are whether or not the last strand is
involved in the braiding. Figure \ref{fig:a13} shows an example of
the two generators of $BF$ corresponding to a generator $A_{ij}$.

The proof of the following proposition is analogous to the proof
of Proposition \ref{bvgens}.

\begin{prop}\label{gens} The elements $x_i$, for $i\ge 0$, $\alpha_{ij}$, for $1\le i<j$, and
$\beta_{ij}$, for $1\le i<j$, form a set of generators of $BF$.
\end{prop}

In the next theorem we will give a presentation for the group
$BF$. The relators are going to be divided into four families. The
family (A) is obtained from the relators of $F$. The family (B) is
obtained from the presentation of the pure braid group:
\begin{itemize}
\item $A_{rs}^{-1}A_{ij}A_{rs}=A_{ij},\text{ if }1\le r<s<i<j\le
n\text{ or }1\le i<r<s<j\le n$ \item
$A_{rs}^{-1}A_{ij}A_{rs}=A_{rj}A_{ij}A_{rj}^{-1},\text{ if }1\le
r<s=i<j\le n$ \item
$A_{rs}^{-1}A_{ij}A_{rs}=(A_{ij}A_{sj})A_{ij}(A_{ij}A_{sj})^{-1},\text{
if }1\le r=i<s<j\le n$ \item $A_{rs}^{-1}A_{ij}A_{rs}=
(A_{rj}A_{sj}A_{rj}^{-1}A_{sj}^{-1})A_{ij}(A_{rj}A_{sj}A_{rj}^{-1}A_{sj}^{-1})^{-1}
,\text{ if }1\le r<i<s<j\le n$
\end{itemize}
obtained from \cite{Hansen}. The family (D) reflects the
interactions between generators of $F$ and pure braids. In
\cite{brinbv1}, Brin constructs these relators using the structure
of Zappa-Sz{\'e}p product of the monoid associated to
$\widehat{BV}$. This construction is not possible here because the
presentations for the groups $P_\infty$ are not monoid
presentations. In fact, the monoid of pure positive braids is not
finitely generated as shown by Burillo, Gutierrez, Krsti{\'c} and
Nitecki. \cite{bgkn}.

The family (C) of relators is given by the special way that the
pure braid groups are embedded into each other inside $BF$. To
embed $P_n$ into $P_{n+1}$ we split the last strand in two. If the
last strand is not braided, this does not affect the element, but
if the last strand takes part in the actual braiding, then these
elements in $P_n$ change when embedded in $P_{n+1}$. When a
generator $\beta_{ij}$ has its last strand split, now the $i$-th
strand wraps around two strands (the $j$-th and the $(j+1)$-th),
and the element is now a product of two generators.

\begin{thm}
The group $BF$ admits a presentation with generators:
\begin{itemize}
\item $x_i$, for $i\ge 0$, \item $\alpha_{ij}$, for $1\le i<j$,
\item $\beta_{ij}$, for $1\le i<j$
\end{itemize}
and relators:
\begin{itemize}
\item[(A)] $x_jx_i=x_ix_{j+1},\text{ if } i<j$ \item[(B1)]
$\alpha_{rs}^{-1}\alpha_{ij}\alpha_{rs}=\alpha_{ij},\text{ if }1\le
r<s<i<j \text{ or }1\le i<r<s<j$ \item[(B2)]
$\alpha_{rs}^{-1}\alpha_{ij}\alpha_{rs}=\alpha_{rj}\alpha_{ij}\alpha_{rj}^{-1}
,\text{ if }1\le r<s=i<j$ \item[(B3)]
$\alpha_{rs}^{-1}\alpha_{ij}\alpha_{rs}=(\alpha_{ij}\alpha_{sj})\alpha_{ij}
(\alpha_{ij}\alpha_{sj})^{-1},\text{ if }1\le r=i<s<j$ \item[(B4)]
$\alpha_{rs}^{-1}\alpha_{ij}\alpha_{rs}=
(\alpha_{rj}\alpha_{sj}\alpha_{rj}^{-1}\alpha_{sj}^{-1})\alpha_{ij}
(\alpha_{rj}\alpha_{sj}\alpha_{rj}^{-1}\alpha_{sj}^{-1})^{-1}
,\text{ if }1\le r<i<s<j$ \item[(B5)]
$\alpha_{rs}^{-1}\beta_{ij}\alpha_{rs}=\beta_{ij},\text{ if }1\le
r<s<i<j \text{ or }1\le i<r<s<j$ \item[(B6)]
$\alpha_{rs}^{-1}\beta_{ij}\alpha_{rs}=\beta_{rj}\beta_{ij}\beta_{rj}^{-1}
,\text{ if }1\le r<s=i<j$ \item[(B7)]
$\alpha_{rs}^{-1}\beta_{ij}\alpha_{rs}=(\beta_{ij}\beta_{sj})
\beta_{ij}(\beta_{ij}\beta_{sj})^{-1},\text{ if }1\le r=i<s<j$
\item[(B8)] $\alpha_{rs}^{-1}\beta_{ij}\alpha_{rs}=
(\beta_{rj}\beta_{sj}\beta_{rj}^{-1}\beta_{sj}^{-1})\beta_{ij}
(\beta_{rn}\beta_{sj}\beta_{rj}^{-1}\beta_{sj}^{-1})^{-1} ,\text{
if }1\le r<i<s<j$ \item[(C)]
$\beta_{ij}=\beta_{i,j+1}\alpha_{ij},\text{ if }i<j$ \item[(D1)]
$\alpha_{ij}x_k=x_k\alpha_{i+1,j+1},\text{ if }k<i-1$ \item[(D2)]
$\alpha_{ij}x_k=x_k\alpha_{i+1,j+1}\alpha_{i,j+1},\text{ if
}k=i-1$ \item[(D3)] $\alpha_{ij}x_k=x_k\alpha_{i,j+1},\text{ if
}i-1<k<j-1$ \item[(D4)]
$\alpha_{ij}x_k=x_k\alpha_{i,j+1}\alpha_{ij},\text{ if }k=j-1$
\item[(D5)] $\alpha_{ij}x_k=x_k\alpha_{ij},\text{ if }k>j-1$
\item[(D6)] $\beta_{ij}x_k=x_k\beta_{i+1,j+1},\text{ if }k<i-1$
\item[(D7)] $\beta_{ij}x_k=x_k\beta_{i+1,j+1}\beta_{i,j+1},\text{
if }k=i-1$ \item[(D8)] $\beta_{ij}x_k=x_k\beta_{i,j+1},\text{ if
}i-1<k<j-1$ \item[(D9)] $\beta_{ij}x_k=x_k\beta_{ij},\text{ if
}k\ge j-1$
\end{itemize}
\end{thm}

\begin{figure}
\includegraphics[width=5in]{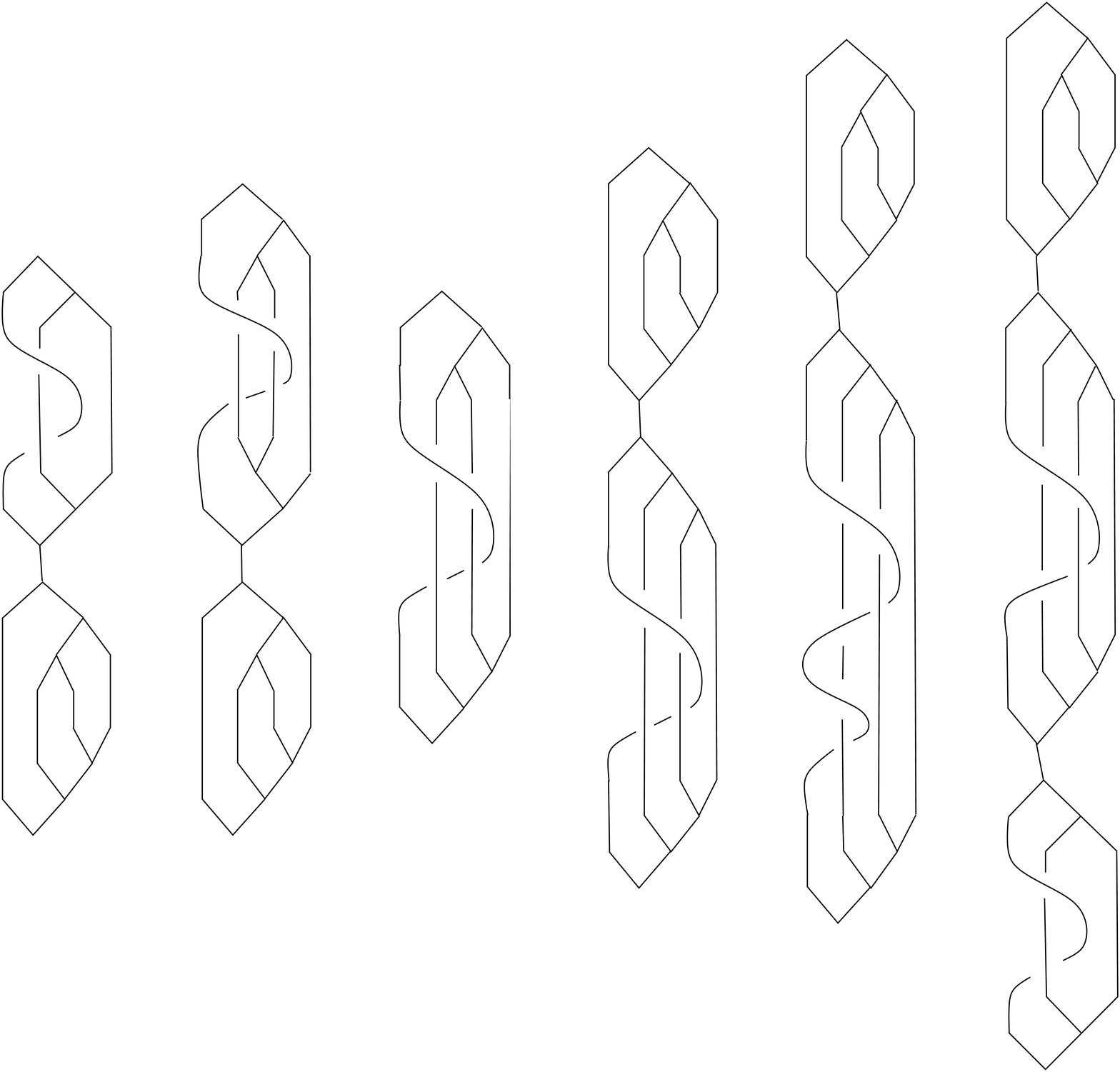}\\
\caption{The process of checking the relator (D4). It starts on
the left with $\alpha_{12}x_1$ and the following steps are the
multiplication processes which transform it into
$x_1\alpha_{13}\alpha_{12}.$\label{fig:relator}}
\end{figure}

{\it Proof.\/} As in the proof of Theorem \ref{bvinfpres}, we
consider the algebraic and the geometric group, establish a
homomorphism between them, which is well-defined and surjective.
For instance, one needs to check geometrically the relators to see
it is well defined. See Figure \ref{fig:relator} for an example.
After this, it only remains to check the injectivity of the map.

To prove that the homomorphism is injective, we consider an
element of $G$, given as a word $w(x_i,\alpha_{ij},\beta_{ij})$,
and imagine that it is mapped to the identity in $G$. We need to
prove that it is consequence of the relators listed above.

As a first step, we can see that the relators (D1) to (D9) and (A)
can be used to transform any such word into the product of three
words,
$$
w_1(x_i) \,w_2(\alpha_{ij},\beta_{ij})\, w_3(x_i)^{-1}.
$$
We can  arrange this in such a way that the words $w_1$ and $w_3$
contain only generators $x_i$ and not their inverses.

This special expression is then particularly useful for studying
the element, because it corresponds easily to the representative
$(T_-,p,T_+)$. We will use now the following lemma whose proof is
analogous to the proof of Lemma \ref{blockid}:

\begin{lem}
A triple $(T_-,p,T_+)$ in $BF$ represents the identity element if
and only if the braid $p$ is the trivial braid and the element
$(T_-,T_+)$ represents the identity in $F$.
\end{lem}

With this lemma, we can assume now that the two words $w_1
w_3^{-1}$ and $w_2$ map to the identity, and we must prove that
they are consequence of the relators. The word $w_1 w_3^{-1}$ lies
in the subgroup isomorphic to $F$. So if it is the identity, it is
consequence of the relators (A).

The word $w_2$ is a product of some $\alpha$ and $\beta$
generators. We would like to consider this word inside some $P_n$,
for a fixed $n$. The image in $BF$ of $\beta_{ij}$ has $j$
strands, and the image of $\alpha_{ij}$ has $j+1$ strands. So the
appropriate $n$ to use is the maximum of the following set:
$$
\{j\,|\,\beta_{ij}\text{ appears in }
w_2\}\cup\{j+1\,|\,\alpha_{ij}\text{ appears in }w_2\}.
$$
If we have a $\beta_{ij}$ with $j<n$, we use the relators (C) to
increase the second index  of that $\beta_{ij}$ to $n$. This way,
the only generators involved are $\alpha_{ij}$, for $1<i<j<n$, and
$\beta_{in}$, for $1<i<n$, which generate a copy of $P_n$ inside
$BF$. So, we can use the relators (C) to have our word expressed
in this small set of generators and assume that it is a word in
$P_n$. So if the word is the identity, it is consequence of the
relators of $P_n$. But these relators correspond to the relators
(B1) to (B8).\qed

\section{A finite presentation for $BF$}

As is common in the groups of the Thompson family, the infinite
presentations are interesting and useful because of their symmetry
and associated normal forms, but often it turns out that there are
finite presentations  from which the infinitely many generators
and relations can be constructed and deduced. In this section, we
construct a finite presentation for $BF$.

Thompson's group $F$ admits a finite presentation which is merely
the first two generators $x_0$ and $x_1$ and the first two
non-trivial relations. We construct $x_n$ from $x_0$ and $x_1$ as
$x_n = x_0^{-n+1} x_1 x_0^{n-1}$ and from the two first two
non-trivial relations we can deduce all of the relations in (A)
above.  This is the first building block for our finite
presentation.

In a similar way we can construct all generators $\alpha_{ij}$
from a few ones. The generators needed are $\alpha_{12}$,
$\alpha_{13}$, $\alpha_{23}$, $\alpha_{24}$. The idea is that
conjugating a braid with $x_k$ has the effect of splitting the
$(k-1)$ strand, so from these four generators and the generators
for $F$, we can construct any generator $\alpha_{ij}$ by the
process of splitting as many strands as necessary to produce the
strands before the $i$-th and between the $i$-th and $j$-th. This
process is as follows:
\begin{itemize}
\item Given $\alpha_{i+1,j+1}$, with $j\ge i+3$, we use the
relators (D3) to decrease the distance between $i$ and $j$ until
2:
$$
\alpha_{i,j+1}=x_{j-2}^{-1}\alpha_{ij}x_{j-2}.
$$
If $i=1$ then this process brings any generator $\alpha_{1j}$ down
to $\alpha_{13}$. For any other value of $i$ it reduces to
$\alpha_{i,i+2}$. \item We reduce the generators
$\alpha_{i+1,i+3}$ with $i\ge 2$ to $\alpha_{24}$ with relators of
type (D1):
$$
\alpha_{i+1,i+3}=x_{i-2}^{-1}\alpha_{i,i+2}x_{i-2}.
$$
\item And finally, we reduce generators of type $\alpha_{i+1,i+2}$
to $\alpha_{23}$ by again using (D1):
$$
\alpha_{i+1,i+2}=x_{i-2}^{-1}\alpha_{i,i+1}x_{i-2}.
$$
\end{itemize}

The generators $\beta$ are constructed in exactly the same way,
where we replace $\alpha$ by $\beta$, and with the same
constraints on indices.

To see which relators to include in the finite presentation, we
see which can be used to get the full families.  We consider the
relators (D) first and we will use those to help with the other
families.

In each of the relators of (D1), there are three strands which are
important: the strands labelled $i$, $j$ and $k$. The idea is that
between those, we only need to have one strand, because by
splitting it, we can get to any number of strands in that
position. Thus, we can get all of the relators (D1) from the
following
\begin{itemize}
\item[(d1.1)] $\alpha_{34}=x_0^{-1}\alpha_{23}x_0$ \item[(d1.2)]
$\alpha_{35}=x_0^{-1}\alpha_{24}x_0$ \item[(d1.3)]
$\alpha_{45}=x_0^{-1}\alpha_{34}x_0$ \item[(d1.4)]
$\alpha_{46}=x_0^{-1}\alpha_{35}x_0$ \item[(d1.5)]
$\alpha_{45}=x_1^{-1}\alpha_{34}x_1$ \item[(d1.6)]
$\alpha_{46}=x_1^{-1}\alpha_{35}x_1$ \item[(d1.7)]
$\alpha_{56}=x_1^{-1}\alpha_{45}x_1$ \item[(d1.8)]
$\alpha_{57}=x_1^{-1}\alpha_{46}x_1$
\end{itemize}
Every relator of the type (D1) is a consequence of the definitions
above and of these eight relators. As an example, we will show the
relator $\alpha_{i,i+1}x_0=x_0\alpha_{i+1,i+2}$. If $i=2$, the
relator is (D1.1) above. If $i>2$, then we use the definitions:
$$
\alpha_{i+1,i+2}=x_{i-2}^{-1}x_{i-3}^{-1}\ldots
x_2^{-1}\alpha_{45}x_2\ldots x_{i-3}x_{i-2}
$$
And using (d1.3) we get to
$$
x_{i-2}^{-1}x_{i-3}^{-1}\ldots x_2^{-1}\ (x_0^{-1}\alpha_{34}x_0)\
x_2\ldots x_{i-3}x_{i-2}
$$
which is
$$
x_0^{-1}x_{i-3}^{-1}x_{i-4}^{-1}\ldots
x_1^{-1}\alpha_{34}x_1\ldots x_{i-4}x_{i-3}x_0
$$
finally equal to
$$
x_0^{-1}\alpha_{i,i+1}x_0.
$$

All the other relators of type (D) are very similar to this case,
and we leave the details to  the reader as they are
straightforward but tedious.

The families of relators (D1) and (D3) are especially important,
because they are used to split or combine adjacent strands which
are not involved in the braiding. Any two adjacent strands which
are not braided can be joined using a relator from one of these
families. This is useful for the families (B1) to (B9).

\begin{figure}
\includegraphics[width=1.5in]{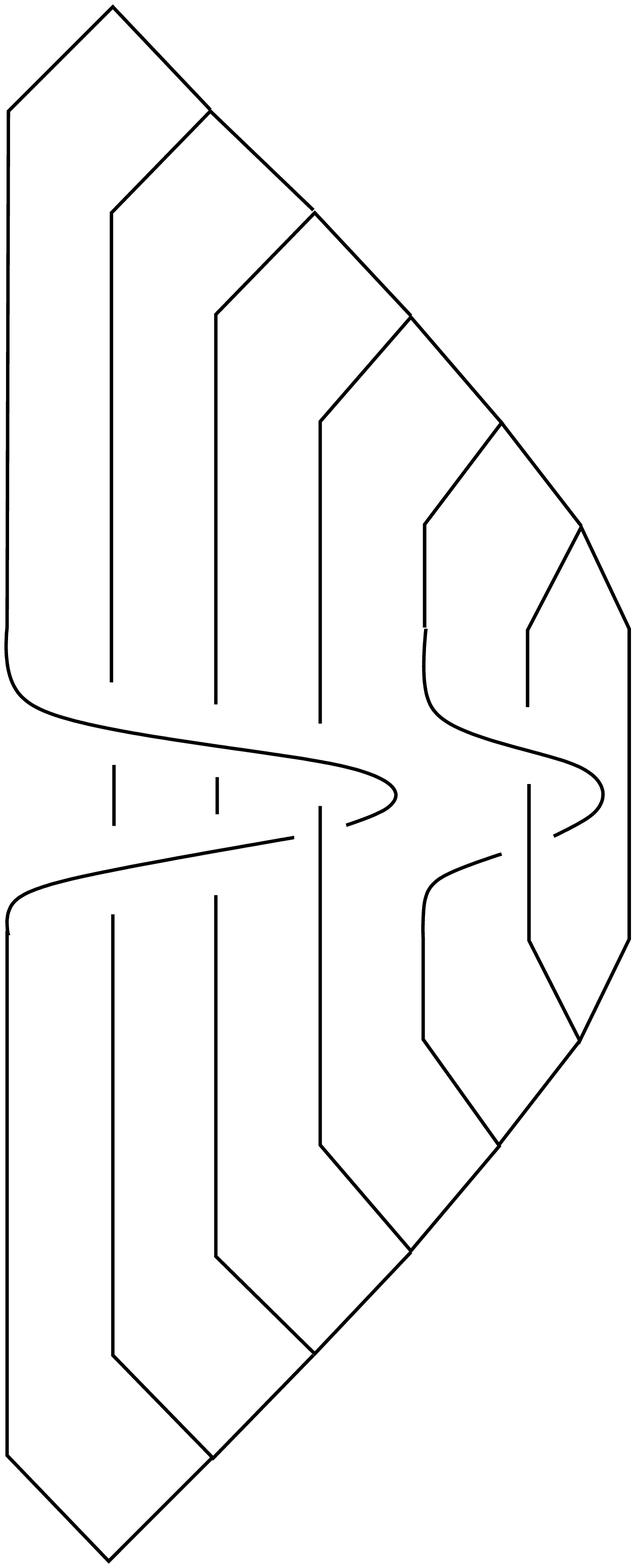}\\
\caption{The element involved in the relator
$\alpha_{14}\alpha_{56}=\alpha_{56}\alpha_{14}$. The commutativity
is apparent. \label{fig:a1456}}
\end{figure}

For instance,  the relators (B1) show that two braids
$\alpha_{ij}$ and $\alpha_{rs}$ commute if $1\le r<s<i<j$. This
relation only involves the strands $r$, $s$, $i$ and $j$. If there
are strands in between, they are uninvolved in the braiding.  If
there is more than one  strand, we can apply the conjugating
relations (D1) or (D3) to bring the relators down to a simple one
where there is just one strand between. For instance, we show the
relator
$$
\alpha_{14}\alpha_{56}=\alpha_{56}\alpha_{14}
$$
in Figure \ref{fig:a1456}.
We see that the two strands between the first and the fourth are
straight. It is clear they can be obtained by splitting a single
strand from an analogous relator whose braided strands are 1, 3, 4
and 5. A conjugation then by $x_1$ brings it down, according to
the definitions above. Now we have that
$$
\alpha_{14}=x_1^{-1}\alpha_{13}x_1,
$$
which is of the type (D3), and
$$
\alpha_{56}=x_1^{-1}\alpha_{45}x_1,
$$
which is of the type (D1). So the relator is a consequence of the
relator
$$
\alpha_{13}\alpha_{45}=\alpha_{45}\alpha_{13},
$$
using only (D1) and (D3) relators. In this way we see that the
only relators that we need to construct all the relators in (B1)
are those that have either zero or one strand at the beginning, or
between the $i$, $j$, $r$, and $s$-th strands. These are:
\begin{itemize}
\item $\alpha_{12}\alpha_{34}=\alpha_{34}\alpha_{12}$ \item
$\alpha_{12}\alpha_{35}=\alpha_{35}\alpha_{12}$ \item
$\alpha_{12}\alpha_{45}=\alpha_{45}\alpha_{12}$ \item
$\alpha_{12}\alpha_{46}=\alpha_{46}\alpha_{12}$ \item
$\alpha_{13}\alpha_{45}=\alpha_{45}\alpha_{13}$ \item
$\alpha_{13}\alpha_{46}=\alpha_{46}\alpha_{13}$ \item
$\alpha_{13}\alpha_{56}=\alpha_{56}\alpha_{13}$ \item
$\alpha_{13}\alpha_{57}=\alpha_{57}\alpha_{13}$ \item
$\alpha_{23}\alpha_{45}=\alpha_{45}\alpha_{23}$ \item
$\alpha_{23}\alpha_{46}=\alpha_{46}\alpha_{23}$ \item
$\alpha_{23}\alpha_{56}=\alpha_{56}\alpha_{23}$ \item
$\alpha_{23}\alpha_{57}=\alpha_{57}\alpha_{23}$ \item
$\alpha_{24}\alpha_{56}=\alpha_{56}\alpha_{24}$ \item
$\alpha_{24}\alpha_{57}=\alpha_{57}\alpha_{24}$ \item
$\alpha_{24}\alpha_{67}=\alpha_{67}\alpha_{24}$ \item
$\alpha_{24}\alpha_{68}=\alpha_{68}\alpha_{24}$
\end{itemize}
The other families, including the ones which involve $\beta$
generators, are completely analogous.

We are left only with the family (C). The relators are all of the
type $\beta_{ij}=\beta_{i,j+1}\alpha_{ij}$, with $i<j$. As before,
it is clear that if $i\ge 3$, we can conjugate by $x_0$ to get a
relator with lower indices. So only ones with $i=1$ or $i=2$ are
needed. Again we see that all are consequences of a fundamental
finite set of relations:
\begin{itemize}
\item $\beta_{12}=\beta_{13}\alpha_{12}$ \item
$\beta_{13}=\beta_{14}\alpha_{13}$ \item
$\beta_{14}=\beta_{15}\alpha_{14}$ \item
$\beta_{23}=\beta_{24}\alpha_{23}$ \item
$\beta_{24}=\beta_{25}\alpha_{24}$ \item
$\beta_{25}=\beta_{26}\alpha_{25}$
\end{itemize}
For example, we see that
$$
\beta_{1j}=x_{j-3}^{-1}\beta_{1,j-1}x_{j-3}=
x_{j-3}^{-1}\beta_{1j}\alpha_{1,j-1}x_{j-3}=\beta_{1,j+1}\alpha_{1j}
$$
if $j\ge 5$. And the cases for $i=2$ are exactly similar.

Hence, we have proved the following theorem:

\begin{thm} \label{bffp}The group $BF$ is finitely presented\end{thm}

The finite presentation for $BF$ is the following:

Generators: $x_0$, $x_1$, $\alpha_{12}$, $\alpha_{13}$,
$\alpha_{23}$, $\alpha_{24}$, $\beta_{12}$, $\beta_{13}$,
$\beta_{23}$, $\beta_{24}$

Relators:
\begin{itemize}
\item[(A)] $x_jx_i=x_ix_{j+1},\text{ for } (i,j)=(1,2),(1,3)$
\item[(B1)] $\alpha_{rs}^{-1}\alpha_{ij}\alpha_{rs}=\alpha_{ij}$,
for
$$
(r,s,i,j)=\begin{cases}\begin{matrix}
(1,2,3,4)&(1,2,3,5)&(1,2,4,5)&(1,2,4,6)\\
(1,3,4,5)&(1,3,4,6)&(1,3,5,6)&(1,3,5,7)\\
(2,3,4,5)&(2,3,4,6)&(2,3,5,6)&(2,3,5,7)\\
(2,4,5,6)&(2,4,5,7)&(2,4,6,7)&(2,4,6,8)\\
(2,3,1,4)&(2,3,1,5)&(2,4,1,5)&(2,4,1,6)\\
(3,4,1,5)&(3,4,1,6)&(3,5,1,6)&(3,5,1,7)\\
(3,4,2,5)&(3,4,2,6)&(3,5,2,6)&(3,5,2,7)\\
(4,5,2,6)&(4,5,2,7)&(4,6,2,7)&(4,6,2,8)\\
\end{matrix}\end{cases}
$$
\item[(B2)]
$\alpha_{rs}^{-1}\alpha_{ij}\alpha_{rs}=\alpha_{rj}\alpha_{ij}\alpha_{rj}^{-1}
$, for
$$
(r,s,i,j)=\begin{cases}\begin{matrix}
(1,2,2,3)&(1,2,2,4)&(1,3,3,4)&(1,3,3,5)\\
(2,3,3,4)&(2,3,3,5)&(2,4,4,5)&(2,4,4,6)
\end{matrix}\end{cases}
$$
\item[(B3)]
$\alpha_{rs}^{-1}\alpha_{ij}\alpha_{rs}=(\alpha_{ij}\alpha_{sj})\alpha_{ij}
(\alpha_{ij}\alpha_{sj})^{-1}$, for
$$
(r,s,i,j)=\begin{cases}\begin{matrix}
(1,2,1,3)&(1,2,1,4)&(1,3,1,4)&(1,3,1,5)\\
(2,3,2,4)&(2,3,2,5)&(2,4,2,5)&(2,4,2,6)
\end{matrix}\end{cases}
$$
\item[(B4)] $\alpha_{rs}^{-1}\alpha_{ij}\alpha_{rs}=
(\alpha_{rj}\alpha_{sj}\alpha_{rj}^{-1}\alpha_{sj}^{-1})\alpha_{ij}
(\alpha_{rj}\alpha_{sj}\alpha_{rj}^{-1}\alpha_{sj}^{-1})^{-1}$,
for
$$
(r,s,i,j)=\begin{cases}\begin{matrix}
(1,3,2,4)&(1,3,2,5)&(1,4,2,5)&(1,4,2,6)\\
(1,4,3,5)&(1,4,3,6)&(1,5,3,6)&(1,5,3,7)\\
(2,4,3,5)&(2,4,3,6)&(2,5,3,6)&(2,5,3,7)\\
(2,5,4,6)&(2,5,4,7)&(2,6,4,7)&(2,6,4,8)\\
\end{matrix}\end{cases}
$$
\item[(B5)] $\alpha_{rs}^{-1}\beta_{ij}\alpha_{rs}=\beta_{ij}$,
for
$$
(r,s,i,j)=\begin{cases}\begin{matrix}
(1,2,3,4)&(1,2,3,5)&(1,2,4,5)&(1,2,4,6)\\
(1,3,4,5)&(1,3,4,6)&(1,3,5,6)&(1,3,5,7)\\
(2,3,4,5)&(2,3,4,6)&(2,3,5,6)&(2,3,5,7)\\
(2,4,5,6)&(2,4,5,7)&(2,4,6,7)&(2,4,6,8)\\
(2,3,1,4)&(2,3,1,5)&(2,4,1,5)&(2,4,1,6)\\
(3,4,1,5)&(3,4,1,6)&(3,5,1,6)&(3,5,1,7)\\
(3,4,2,5)&(3,4,2,6)&(3,5,2,6)&(3,5,2,7)\\
(4,5,2,6)&(4,5,2,7)&(4,6,2,7)&(4,6,2,8)\\
\end{matrix}\end{cases}
$$
\item[(B6)]
$\alpha_{rs}^{-1}\beta_{ij}\alpha_{rs}=\beta_{rj}\beta_{ij}\beta_{rj}^{-1}$
, for
$$
(r,s,i,j)=\begin{cases}\begin{matrix}
(1,2,2,3)&(1,2,2,4)&(1,3,3,4)&(1,3,3,5)\\
(2,3,3,4)&(2,3,3,5)&(2,4,4,5)&(2,4,4,6)
\end{matrix}\end{cases}
$$
\item[(B7)]
$\alpha_{rs}^{-1}\beta_{ij}\alpha_{rs}=(\beta_{ij}\beta_{sj})
\beta_{ij}(\beta_{ij}\beta_{sj})^{-1}$, for
$$
(r,s,i,j)=\begin{cases}\begin{matrix}
(1,2,1,3)&(1,2,1,4)&(1,3,1,4)&(1,3,1,5)\\
(2,3,2,4)&(2,3,2,5)&(2,4,2,5)&(2,4,2,6)
\end{matrix}\end{cases}
$$
\item[(B8)] $\alpha_{rs}^{-1}\beta_{ij}\alpha_{rs}=
(\beta_{rj}\beta_{sj}\beta_{rj}^{-1}\beta_{sj}^{-1})\beta_{ij}
(\beta_{rn}\beta_{sj}\beta_{rj}^{-1}\beta_{sj}^{-1})^{-1}$, for
$$
(r,s,i,j)=\begin{cases}\begin{matrix}
(1,3,2,4)&(1,3,2,5)&(1,4,2,5)&(1,4,2,6)\\
(1,4,3,5)&(1,4,3,6)&(1,5,3,6)&(1,5,3,7)\\
(2,4,3,5)&(2,4,3,6)&(2,5,3,6)&(2,5,3,7)\\
(2,5,4,6)&(2,5,4,7)&(2,6,4,7)&(2,6,4,8)\\
\end{matrix}\end{cases}
$$
\item[(C)] $\beta_{ij}=\beta_{i,j+1}\alpha_{ij}$, for
$$
(i,j)=(1,2),(1,3),(1,4),(2,3),(2,4),(2,5)
$$
\item[(D1)] $\alpha_{ij}x_k=x_k\alpha_{i+1,j+1}$, for
$$
(i,j,k)=\begin{cases}\begin{matrix}
(2,3,0)&(2,4,0)&(3,4,0)&(3,5,0)\\
(3,4,1)&(3,5,1)&(4,5,1)&(4,6,1)
\end{matrix}\end{cases}
$$
\item[(D2)] $\alpha_{ij}x_k=x_k\alpha_{i+1,j+1}\alpha_{i,j+1}$,
for
$$
(i,j,k)=(1,2,0),(1,3,0),(2,3,1),(2,4,1)
$$
\item[(D3)] $\alpha_{ij}x_k=x_k\alpha_{i,j+1}$, for
$$
(i,j,k)=\begin{cases}\begin{matrix}
(1,3,1)&(1,4,1)&(1,4,2)&(1,5,2)\\
(2,4,2)&(2,5,2)&(2,5,3)&(2,6,3)
\end{matrix}\end{cases}
$$
\item[(D4)] $\alpha_{ij}x_k=x_k\alpha_{i,j+1}\alpha_{ij}$, for
$$
(i,j,k)=(1,2,1),(1,3,2),(2,3,2),(2,4,3)
$$
\item[(D5)] $\alpha_{ij}x_k=x_k\alpha_{ij}$, for
$$
(i,j,k)=\begin{cases}\begin{matrix}
(1,2,2)&(1,2,3)&(1,3,3)&(1,3,4)\\
(2,3,3)&(2,3,4)&(2,4,4)&(2,4,5)
\end{matrix}\end{cases}
$$
\item[(D6)] $\beta_{ij}x_k=x_k\beta_{i+1,j+1}$, for
$$
(i,j,k)=\begin{cases}\begin{matrix}
(2,3,0)&(2,4,0)&(3,4,0)&(3,5,0)\\
(3,4,1)&(3,5,1)&(4,5,1)&(4,6,1)
\end{matrix}\end{cases}
$$
\item[(D7)] $\beta_{ij}x_k=x_k\beta_{i+1,j+1}\beta_{i,j+1}$, for
$$
(i,j,k)=(1,2,0),(1,3,0),(2,3,1),(2,4,1)
$$
\item[(D8)] $\beta_{ij}x_k=x_k\beta_{i,j+1}$, for
$$
(i,j,k)=\begin{cases}\begin{matrix}
(1,3,1)&(1,4,1)&(1,4,2)&(1,5,2)\\
(2,4,2)&(2,5,2)&(2,5,3)&(2,6,3)
\end{matrix}\end{cases}
$$
\item[(D9)] $\beta_{ij}x_k=x_k\beta_{ij}$, for
$$
(i,j,k)=(1,2,1),(1,3,2),(2,3,2),(2,4,3)
$$
\end{itemize}

This gives a total of 10 generators and 192 relators.

\section{The braided Thompson group $\widehat{BF}$}

In the previous section we constructed a presentation of $BF$.
Brin \cite{brinbv1} described both $BV$ and $\widehat{BV}$.  Brin
describes the group  $\widehat{BV}$ via a Zappa-Sz\'ep product of
$F$ and $B_\infty$, and describes $BV$ as a subgroup of
$\widehat{BV}$.  The group $\widehat{BV}$ is the group of braided
forest diagrams, where all but finitely many of the forests are
trivial, and the group $BV$ is the subgroup of $\widehat{BV}$
where all of the trees in the forest pairs are trivial except the
first pair, which have the same number of leaves. Not only is $BV$
a subgroup of $\widehat{BV}$, but also $\widehat{BV}$ is a
subgroup of $BV$, as described by Brin \cite{brinbv1}.  We take
the standard identification of the real line with the unit
interval which is compatible with the relevant dyadic subdivisions
which sends the interval $[i,i+1]$ of ${\bf R}$ with the interval
$[1-2^{-i}, 1-2^{-i-1}]$ and then we see that $\widehat{BV}$ is
the subgroup of $BV$ in which the last strand is not braided with
any other strands.  Similarly, we have the group $\widehat{BF}$
which can either be regarded as the supergroup  of $BF$ of pure
braided forest diagrams, or as a subgroup of $BF$ where the
braiding does not involve the last strand.

Here, we easily describe the subgroup $\widehat{BF}$ of $BF$ by
omitting the generators and relations from $BF$ which involve
braiding the last strand.  So we obtain presentations for
$\widehat{BF}$ which are sub-presentations of the  infinite and
finite presentations for $BF$ given in the earlier sections above.

\begin{prop}\label{gensbfhat} The elements $x_i$, for $i\ge 0$, $\alpha_{ij}$, for $1\le i<j$
form a set of generators of $\widehat{BF}$.
\end{prop}

{\it Proof.\/} The proof is similar to the case for $BF$.  Here we
note that  $\widehat{BF}$ is exactly the subgroup where the last
strand is not braided with any previous strands, and by omitting
the $\beta$ generators we guarantee that the last strand is not
braided.  An argument similar to the earlier one for $BF$ shows
that these generate $\widehat{BF}$. \qed

To find relators for a presentation of $\widehat{BF}$ we use the
same sets of generators and relators as for $BF$, deleting the
generators in the $\beta$ family and deleting all relations which
include any of the $\beta_{ij}$.

We thus obtain the following:

\begin{thm}
The group $\widehat{BF}$ admits a presentation with generators:
\begin{itemize}
\item $x_i$, for $i\ge 0$, \item $\alpha_{ij}$, for $1\le i<j$,
\end{itemize}
and relators:
\begin{itemize}
\item[(A)] $x_jx_i=x_ix_{j+1},\text{ if } i<j$ \item[(B1)]
$\alpha_{rs}^{-1}\alpha_{ij}\alpha_{rs}=\alpha_{ij},\text{ if
}1\le r<s<i<j \text{ or }1\le i<r<s<j$ \item[(B2)]
$\alpha_{rs}^{-1}\alpha_{ij}\alpha_{rs}=\alpha_{rj}\alpha_{ij}\alpha_{rj}^{-1}
,\text{ if }1\le r<s=i<j$ \item[(B3)]
$\alpha_{rs}^{-1}\alpha_{ij}\alpha_{rs}=(\alpha_{ij}\alpha_{sj})\alpha_{ij}
(\alpha_{ij}\alpha_{sj})^{-1},\text{ if }1\le r=i<s<j$ \item[(B4)]
$\alpha_{rs}^{-1}\alpha_{ij}\alpha_{rs}=
(\alpha_{rj}\alpha_{sj}\alpha_{rj}^{-1}\alpha_{sj}^{-1})\alpha_{ij}
(\alpha_{rj}\alpha_{sj}\alpha_{rj}^{-1}\alpha_{sj}^{-1})^{-1}
,\text{ if }1\le r<i<s<j$
\item[(D1)] $\alpha_{ij}x_k=x_k\alpha_{i+1,j+1},\text{ if }k<i-1$
\item[(D2)]
$\alpha_{ij}x_k=x_k\alpha_{i+1,j+1}\alpha_{i,j+1},\text{ if
}k=i-1$ \item[(D3)] $\alpha_{ij}x_k=x_k\alpha_{i,j+1},\text{ if
}i-1<k<j-1$ \item[(D4)]
$\alpha_{ij}x_k=x_k\alpha_{i,j+1}\alpha_{ij},\text{ if }k=j-1$
\item[(D5)] $\alpha_{ij}x_k=x_k\alpha_{ij},\text{ if }k>j-1$
\end{itemize}
\end{thm}

{\it Proof.\/}Using the interpretation of the geometric group
$\widehat{BF}$ as the subgroup of $BF$ where braiding never
involves the last strand, the same proof used in the section 5 to
establish the presentation for $BF$ goes through in this
situation. The only difference is that since we have no $\beta$
generators, once we rearrange the word to have all $\alpha$
generators in the middle, they are already generators for one copy
of the pure braids on $n$ strands with the rightmost strand
unbraided, so there is no need for the step using relators of type
(C). \qed

Note that the group $\widehat{BF}$ will also be finitely
presented. The arguments needed to see this are similar to those
for $BF$. The finite presentation for $\widehat{BF}$ can be easily
obtained from the finite presentation for $BF$ by deleting all
generators $\beta_{ij}$ and the relations where they appear.

\bibliographystyle{plain}

\begin{thebibliography}{10}

\bibitem{belkdiss}
James Belk.
\newblock {\em {T}hompson's group ${F}$}.
\newblock PhD thesis, Cornell University, 2004.

\bibitem{brinbv2}
Matthew~G. Brin.
\newblock {The Algebra of Strand Splitting. I. A Braided Version of Thompson's
  Group V}.

\bibitem{brinbv1}
Matthew~G. Brin.
\newblock {The Algebra of Strand Splitting. II. A Presentation for the Braid
  Group on One Strand}.

\bibitem{bgkn}
J.~Burillo, M.~Gutierrez, S.~Krsti{\'c}, and Z.~Nitecki.
\newblock Crossing matrices and {T}hurston's normal form for braids.
\newblock {\em Topology Appl.}, 118(3):293--308, 2002.

\bibitem{bcs}
Jos\'e Burillo, Sean Cleary, and Melanie Stein.
\newblock Metrics and embeddings of generalizations of {T}hompson's group
  ${F}$.
\newblock {\em Trans. Amer. Math. Soc.}, 353(4):1677--1689 (electronic), 2001.

\bibitem{thompt}
Jos\'e Burillo, Sean Cleary, Melanie Stein, and Jennifer Taback.
\newblock Combinatorial and metric properties of {T}hompson's group $t$.
\newblock submitted.

\bibitem{cfp}
J.~W. Cannon, W.~J. Floyd, and W.~R. Parry.
\newblock Introductory notes on {R}ichard {T}hompson's groups.
\newblock {\em Enseign. Math. (2)}, 42(3-4):215--256, 1996.

\bibitem{dehornoy2}
Patrick Dehornoy.
\newblock Geometric presentations for {T}hompson's groups.
\newblock {\em J. Pure Appl. Algebra}, 203(1-3):1--44, 2005.

\bibitem{dehornoy}
Patrick Dehornoy.
\newblock {The group of parenthesized braids}.
\newblock {\em Adv. in Math.}, to appear.

\bibitem{blakegd}
S.~Blake Fordham.
\newblock Minimal length elements of {T}hompson's group {$F$}.
\newblock {\em Geom. Dedicata}, 99:179--220, 2003.

\bibitem{Hansen}
Vagn~Lundsgaard Hansen.
\newblock {\em Braids and coverings: selected topics}, volume~18 of {\em London
  Mathematical Society Student Texts}.
\newblock Cambridge University Press, Cambridge, 1989.
\newblock With appendices by Lars G\ae de and Hugh R. Morton.

\end{thebibliography}

\end{document}